\documentclass[twoside, onecolumn, 10pt]{article}
\usepackage{latexsym}
\usepackage{amssymb}
\usepackage{mathrsfs}
\usepackage{amsmath}
\usepackage{marvosym}

\usepackage{pdfpages}

\newtheorem{theorem}{Theorem}[section]
\newtheorem{lemma}[theorem]{Lemma}
\newtheorem{corollary}[theorem]{Corollary}

\newtheorem{problem}[theorem]{Problem}

\newcommand{\e}{\begin{eqnarray}}
\newcommand{\ee}{\end{eqnarray}}
\newcommand{\en}{\begin{eqnarray*}}
\newcommand{\een}{\end{eqnarray*}}
\newcommand{\nb}{\nonumber}

\def\bc{\begin{center}}
\def\ec{\end{center}}

\begin{document}

\bc
 {\LARGE Analytical formulas for calculating the extremal ranks of the matrix-valued function $A + BXC$
 when the rank of $X$  is fixed}
\ec

\begin{center}
{\large  Yongge Tian}
\end{center}

\begin{center}
{\em {\footnotesize CEMA, Central University of Finance and Economics, Beijing 100081, China
}}
\end{center}

\renewcommand{\thefootnote}{\fnsymbol{footnote}}
\footnotetext{ {\it E-mail Address:} yongge.tian@gmail.com}

\noindent {\small {\bf Abstract.} One of the simplest matrix-valued function with a single variable matrix $X$
is given by $A + BXC$. In this this note, analytical formulas are established for calculating the
 maximal and minimal ranks of $A + BXC$ when the rank of the variable matrix $X$ is fixed by using a simultaneous
 decomposition of $A$,  $B$ and $C$ and some preliminary results. Some applications of the formulas in completing
 partially-specified block matrix with the maximal and minimal ranks  are also given.\\


\medskip

}

\section{Introduction}
\renewcommand{\theequation}{\thesection.\arabic{equation}}
\setcounter{section}{1} \setcounter{equation}{0}

Fixed-rank or low-rank matrix approximation problems are to approximate optimally, with respect to some criteria,
a matrix by one of the same dimension but fixed or smaller rank from a given feasible matrix set.
Assume that $A$ is a matrix to be approximated. Then a conventional statement of general matrix optimization  problems of $A$
from this point of view can be written as
\begin{align}
{\rm minimize}\, \rho(\, A - X \,) \ \  {\rm subject \ to}  \ X  \in {\cal S},
\label{ww11}
\end{align}
where $\rho(\cdot)$ is a certain objective function of decision matrix, which is usually taken as  determinant, trace,
norms, rank, inertia of a matrix, and ${\cal S}$ is a certain feasible matrix set.
A best-known case of (\ref{ww11}) is to minimize the norm $\|\, A - X \,\|^2_{F}$  subject  to  $ X  \in {\cal S}.$
The fixed-rank or low-rank matrix set mentioned above can be written as
\begin{align}
{\cal S} = \{X  \ | \ {\rm rank}(X) = t \} \ \  {\rm or} \ \ {\cal S} = \{ X  \ | \ {\rm rank}(X) \leqslant t \}.
\label{ww12}
\end{align}
The use of low-rank matrix to approximate a given matrix dates back to \cite{EY,Sc}, which now becomes a very active research subject in both optimization theory and applied disciplines.

Although these problems are stated quite clearly in form, it is hard in general to give satisfactory answers
in closed-form to these matrix approximation problems.  In other words, only numerical solutions to these
approximation problems can be obtained. In this note, we assume that the objective function $\rho(\cdot)$ in
(\ref{ww11})  is taken as the rank of  matrix. Then this kind of optimization problems can generally be written as
\begin{align}
& \text{maximize}\ {\rm rank}(\, A - X\,) \ \ \ \ \ \,   {\rm subject \ to}  \ \ \  X \in {\cal S},
\label{ww13}
\\
& \text{minimize}\ {\rm rank}(\, A - X\,) \ \ \ \ \ \   {\rm subject \ to}  \ \ \  X \in {\cal S},
\label{ww14}
\end{align}
respectively. The rank of matrix, as an objective function, is often used when finding feasible matrix $X$ such that resulting
$A - X$ attains its maximal possible rank (is nonsingular when square), or such that $A - X$ attains the minimal rank as
possible (called low-rank matrix completion). This kind of  problems
are usually called the matrix rank-maximization and rank-minimization problems, or  matrix rank completion
problems in the literature. Generally speaking, matrix rank-optimization problems are a class of discontinuous
optimization problems, in which the decision variables are matrices running over certain matrix sets, while
the ranks of the variable matrices are taken as integer-valued objective functions.  In  this case, analytical formulas for
calculating the integer extremum ranks of $A - X$ can hardly be derived by
numerical approximation methods. This fact means that solving methods of matrix rank optimization
 problems are not consistent with any of the ordinary continuous and discrete problems in optimization theory,
so that we cannot apply various common methods of solving continuous
optimization problems, such as the well-known differential and Lagrangian methods, to approach these constrained
optimization problems. Instead, we can only find the exact maximal and minimal ranks
 through pure algebraic operations of matrices.   It has been known that matrix rank-optimization problems
 are NP-hard in general due to the discontinuity and combinational nature of rank of a matrix and the
 algebraic structure of ${\cal S}$. Many new researches were conducted on this kind of matrix rank-optimization problems
 from theory and applied points of  view in the past decades;  see, e.g., \cite{MGC}.
% \texttt{http://perception.csl.uiuc.edu/matrix-rank/home.html} and the references therein.
Because the rank of a matrix can only take finite integers between 0 and the dimensions of the matrix, it is really
expected to establish certain analytical formulas for calculating the maximal and minimal ranks for curiosity.

 In what follows, we assume that $A \in {\mathbb C}^{m \times n}$,  $B \in {\mathbb C}^{m \times p}$ and
 $C \in {\mathbb C}^{q \times n}$ are given matrices,  and the feasible matrix set ${\cal S}$ in (\ref{ww11})
 is taken as
\begin{align}
{\cal S} = \{ - BXC  \ | \ X\in {\mathbb C}^{p\times q} \ {\rm and} \  {\rm rank}(X) = t \}.
\label{ww17}
\end{align}
Then, the difference in (\ref{ww11}) can equivalently be written as the following linear matrix-valued function
\begin{align}
\phi(X) = A + BXC,
\label{11}
\end{align}
which is a map $\phi :  {\mathbb C}^{p \times q} \rightarrow  {\mathbb C}^{m \times n}$. Under such a
formulation, this note aims at solving the following constrained matrix optimization problems:

\begin{problem} \label{T11}
{\rm For the function in (\ref{11}) and two integers $s$ and $t$ with
$0\leqslant s \leqslant t \leqslant \min\{\,p, \ q\,\}$, establish
explicit formulas for calculating the following extremal ranks
\begin{align}
& \text{maximize} \  {\rm rank}(\, A + BXC \,)   \ \ \ \ \ \ \ \  \mbox{s.t.}  \ \  X \in  {\mathbb C}^{p\times q}
\ \ \mbox{and} \ \  {\rm rank}(X) = t,
\label{12}
\\
& \text{minimize} \ {\rm rank}(\, A + BXC \,)  \ \ \ \ \ \ \ \  \mbox{s.t.}  \ \  X \in  {\mathbb C}^{p\times q}
\ \ \mbox{and} \ \   {\rm rank}(X) = t,
\label{13}
\\
& \text{maximize} \ {\rm rank}(\, A + BXC \,)   \ \ \ \ \ \ \ \,  \mbox{s.t.}  \ \  X \in  {\mathbb C}^{p\times q}
\ \ \mbox{and} \ \   s \leqslant  {\rm rank}(X) \leqslant t,
\label{14}
\\
& \text{minimize} \ {\rm rank}(\, A + BXC \,)  \ \ \ \ \ \ \ \  \mbox{s.t.}  \ \  X \in  {\mathbb C}^{p\times q}
\ \ \mbox{and} \ \   s \leqslant  {\rm rank}(X) \leqslant t.
\label{15}
\end{align}
}
\end{problem}

The matrices $X$ satisfying the constraints in (\ref{12})--(\ref{15}) are called the feasible solutions
(i.e., candidates for solutions) to the problems, respectively. They form  certain sets of ${\mathbb C}^{p\times q}$
and it is over these sets that the objective function is to be maximized or minimized. However, these matrix sets are not
 necessarily convex. Motivations for finding the extremal ranks of (\ref{11}) arise from both theoretical and applied
points of view. It is really lucky that we can establish analytical formulas
 for calculating the extremal ranks of matrix-valued functions for some special matrix sets ${\cal S}$
 by using various expansion formulas for ranks of matrices and some tricky matrix operations.
For instance, two well-known seminal formulas in closed-form for calculating the global maximal and minimal
ranks of (\ref{11}) are given by
\begin{align}
&  \max_{X \in \mathbb C^{p \times q}}{\rm rank}(\, A + BXC \,)
 =\min\left\{ {\rm rank}[\,A, \, B\,],\ {\rm rank}\!\left[\!\begin{array}{c} A \\ C
 \end{array}\!\right] \right\},
\label{16}
\\
 & \min_{ X \in \mathbb C^{p \times q}}{\rm rank}(\, A + BXC \,)
 = {\rm rank}[\,A, \, B\,]+ {\rm rank}\!\left[\!\begin{array}{c} A \\ C
 \end{array}\!\right] - {\rm rank}\!\left[\!\begin{array}{cc} A & B \\ C & 0
 \end{array}\!\right].
\label{17}
\end{align}
Because the right-hand sides of (\ref{16}) and (\ref{17}) are calculated only by three block matrices composed by
the three given matrices, a beginner who knows the concept of matrix rank in linear algebra can understand the usefulness
of  (\ref{16}) and (\ref{17}). People can apply (\ref{16}) and (\ref{17}) to characterize many
fundamental behaviors of $A + BXC$, for instance, necessary and sufficient conditions can directly be established
for $A + BXC$ to be nonsingular; for $A + BXC$ to be zero; for the rank of $A + BXC$ to be invariant under
different choice of $X$;  for the row and column spaces of $A + BXC$ to be invariant under different choice of
$X$, respectively, etc. However, these two elementary formulas cannot be proved within the scope of elementary linear
algebra. Some people made essential contributions for the establishments of (\ref{16}) and (\ref{17})
through pure algebraic operations of the given matrices and generalized inverses, as well as simultaneous
matrix decompositions of the given matrices; see, e.g., \cite{DG,LT-lama,Ti-1,TC}. Analytical expressions for
the general expressions of the variable matrices $X$ satisfying (\ref{16}) and (\ref{17}) were also
obtained through generalized inverses and simultaneous matrix decompositions of the given matrices in \cite{LT-lama,TC}.
Eqs. (\ref{16}) and (\ref{17}) are not just two isolated formulas for the maximal and minimal ranks
of matrix-valued functions.
%In the past two decades, analytical formulas were widely established for calculating
%the maximal and minimal ranks of some general linear and quadratic matrix-valued functions and their applications
%were presented; see, e.g., \cite{T-laa02,T-mis,T-lama11}.
Motivated by some recent work on low-rank
matrix approximations,  the present author revisits (\ref{11}) by adding certain rank restrictions on the
variable matrix $X$, and establishes some new and elementary formulas for calculating the maximal and minimal ranks
in (\ref{12})--(\ref{15}), which,  we believe,  can be taken as some standard examples for verifying accuracy
of various algorithms in solving matrix rank-approximation problems.

Throughout this note, ${\mathbb C}^{m\times n}$ stands for the set of all $m\times n$ complex
matrices;  ${\mathbb C}_t^{m\times n}$ stands for the set of all $m\times n$ complex
matrices with ${\rm rank}(X) =t$;  $A^{*}$,  $r(A)$ and ${\mathscr R}(A)$ stand
for the conjugate transpose, rank and range (column space) of a
matrix $A\in {\mathbb C}^{m\times n}$, respectively;  $I_m$ denotes the identity matrix of order $m$;
 $[\, A, \, B\,]$ denotes a row block matrix consisting of $A$ and $B$.

In dealing with  problems in the formats of (\ref{12})--(\ref{15}), people usually
 construct certain canonical forms of the matrix-valued functions through some
 simultaneous decompositions of $A$,  $B$ and $C$, because the ranks of matrices are
 invariant under nonsingular matrix transformations.  In order to establish a
 canonical form of  (\ref{11}), we need the following several known or simple results
 on simultaneous decompositions of matrices and rank formulas for block matrices.

\begin{lemma} [\cite{Zha-1,Zha-2}]\label{T13}
Let $A\in {\mathbb C}^{m\times n},$ $B\in {\mathbb C}^{m\times p}$
and $C\in {\mathbb C}^{q\times n}.$ Then there exist two
 nonsingular matrices $P\in \mathbb C^{m \times m},$ $Q \in \mathbb C^{n \times
 n}$ and two unitary matrices $U\in  \mathbb C^{p \times p},$ $V\in  \mathbb C^{q \times q}$ such that
\begin{equation}
 A=P\Sigma_AQ, \ \   B=P\Sigma_BU,\  C= V\Sigma_CQ ,
\label{112}
\end{equation}
where
\begin{align}
&  \begin{array}{c@{\hspace{-15 pt}}l}
\Sigma_A = \left[
\begin{array}{cccccc}
 I  & 0  & 0  & 0 & 0 & 0
 \\
0  & I & 0  & 0 & 0 & 0
\\
0  & 0  & I  & 0 & 0 & 0
\\
0  & 0  & 0  & S_A & 0 & 0
\\
0  & 0  & 0  & 0 & 0 & 0
\\
0  & 0  & 0  & 0 & 0 & 0
\end{array}
\right] & \ \ \
\begin{array}{l}
j \\ k \\ l \\ u \\ s_2 \\ t_2
\end{array}
\\
{\footnotesize  \ \ \ \ \qquad j \ \quad  k \ \quad   l \ \quad \ u \ \quad  s_1  \ \quad  t_1}
\end{array},
\label{113}
\\
&  \begin{array}{c@{\hspace{-15 pt}}l}
\Sigma_B = \left[
\begin{array}{cccc}
 I  &  \qquad 0 \qquad  \ & 0  & 0
 \\
0  &\qquad  0 \qquad  \ & 0  & 0
\\
0  & \qquad 0\qquad  \ & S_B & 0
\\
0  &\qquad  0 \qquad  \ & 0  & I
\\
0  & \qquad 0 \qquad  \ & 0  & 0
\end{array}
\right] & \ \ \
\begin{array}{l}
j \\ k +l \\ u \\ s_2 \\ t_2
\end{array}
\\
{\footnotesize  \ \ \ \ \qquad j\quad  p-j-r-s_2 \ \ u \quad \ s_2}
\end{array},
\label{114}
\\
& \begin{array}{c@{\hspace{-15 pt}}l}
\Sigma_C = \left[
\begin{array}{ccccc}
 0 \qquad  & 0  & 0  & 0 & 0
 \\
0 \qquad  & I & 0  & 0 & 0
\\
0 \qquad  & 0  & S_C  & 0 & 0
\\
0  \qquad  & 0  & 0  & I & 0
\end{array}
\right] & \ \ \
\begin{array}{l} q - l - u - s_1\\ l \\ u \\ s_1
\end{array}
\\
{\footnotesize  \ \ \ \ \qquad j + k \ \ \quad   l \ \quad \ u \ \ \quad  s_1  \quad  t_1}
\end{array},
\label{115}
\end{align}
$S_A,$ $S_B$ and $S_C$ are diagonal matrices with positive  diagonal entries$,$ and
\begin{align*}
& j =r\!\left[\!\! \begin{array}{cc}  A  \\ C  \end{array}
\!\!\right] + r(B) - r\!\left[\!\! \begin{array}{cc}  A  & B  \\ C  & 0 \end{array}
\!\!\right],
\\
&  k =  r\!\left[\!\! \begin{array}{cc}  A  & B  \\ C  & 0 \end{array}
\!\!\right] - r(B) - r(C),
\\
& l =r[\, A, \, B\,] + r(C) - r\!\left[\!\! \begin{array}{cc}  A  & B  \\ C  & 0 \end{array}
\!\!\right],
\\
& u = r\!\left[\!\! \begin{array}{cc}  A  & B  \\ C  & 0 \end{array}
\!\!\right] + r(A) - r\!\left[\!\! \begin{array}{cc}  A  \\ C  \end{array}
\!\!\right]  - r[\, A, \, B\,],
\\
& s_1 =  r\!\left[\!\! \begin{array}{cc}  A  \\ C  \end{array}
\!\!\right] - r(A),
\\
& s_2 = r[\, A, \, B\,] - r(A),
\\
& t_1 = n - r\!\left[\!\! \begin{array}{cc}  A  \\ C  \end{array}
\!\!\right],
\\
& t_2 = m - r[\, A, \, B\,].
\end{align*}
\end{lemma}

\begin{lemma} \label{T14}
Let $X \in \mathbb C^{m \times n},$  $Y \in \mathbb C^{m \times p}$ and $Z \in \mathbb C^{q \times n}$ be
three variable matrices$,$  and let
\begin{align}
\phi(X, \, Y, \, Z\,) = \left[\!\begin{array}{cc} X & Y \\ Z & 0
\end{array}\!\right]\!.
\label{116}
\end{align}
Then$,$
\begin{align}
& \max_{X\in {\mathbb C}^{m \times n}, \,Y \in \mathbb C^{m \times p},\, Z \in \mathbb C^{q \times n}}
r[\phi(X, \, Y, \, Z\,)]  = \min\{\, m + q,   \ \  n + p,  \ \  m + n \,\},
\label{117}
\\
& \min_{X\in {\mathbb C}^{m \times n}, \,Y \in \mathbb C^{m \times p}, \, Z \in \mathbb C^{q \times n}}
r[\phi(X, \, Y, \, Z\,)] = 0.
\label{118}
\end{align}
Further$,$ for any integer $t$ with $0 \leqslant t \leqslant \min\{\, m + q, \ n + p, \  m + n \,\},$ there
exist $X \in \mathbb C^{m \times n},$ $Y\in \mathbb C^{m \times p}$ and $Z \in \mathbb C^{q \times n}$
such that
\begin{align}
r\left[\!\begin{array}{cc} X & Y \\ Z & 0 \end{array}\!\right] = t.
\label{119}
\end{align}
\end{lemma}

\noindent {\bf Proof.} \
It is obvious that the right-hand side of (\ref{117}) is an upper bound of $r[\phi(X, \, Y, \, Z\,)]$.

(I) Under $m + q \leqslant  \min\{\,  n + p,  \ \  m + n \,\}$ and $m \leqslant p$,   setting
$$
X =0, \ \ \ Y = [\, I_m, \ 0\,], \ \  Z = [\, I_q, \ 0\,]
$$
leads to $r[\phi(X, \, Y, \, Z\,)] = r(Y) + r(Z) = m +q$; under
$m + q \leqslant  \min\{\,  n + p,  \ \  m + n \,\}$ and $m > p$,
setting
$$
[\,X, \, Y \,] = [\, 0, \ I_m\,], \ \  Z = [\, I_q, \ 0\,]
$$
leads to $r[\phi(X, \, Y, \, Z\,)] = m +q$;

(II) under $n + p \leqslant  \min\{\,  m + q,  \ \  m + n \,\}$ and $n \leqslant q$, setting
$$
X =0, \ \ Y = [\, I_n, \ 0\,]^T, \ \  Z = [\, I_p, \ 0\,]^T
$$
leads to $r[\phi(X, \, Y, \, Z\,)] = r(Y) + r(Z) = n +p$; under $n + p \leqslant
\min\{\,  m + q,  \ \  m + n \,\}$ and $n > q$,
setting
$$
\left[\!\begin{array}{cc} X \\ Z \end{array}\!\right]  = \left[\!\begin{array}{cc} 0 \\ I_n
\end{array}\!\right], \ \  Y = \left[\!\begin{array}{cc} I_p \\ 0
\end{array}\!\right]
$$
leads to $r[\phi(X, \, Y, \, Z\,)] = r(Y) + r(Z) = n +p$;

(III) under $m + n \leqslant  \min\{\,  m + q,  \ \  n + p \,\}$,   setting
$$
X =0, \ \ \ Y = [\, I_m, \ 0\,], \ \  Z = [\, I_n, \ 0\,]^T
$$
leads to $r[\phi(X, \, Y, \, Z\,)] = r(Y) + r(Z) =m +n$; establishing  (\ref{117}).

Setting $X = 0$ and $Y = 0$ leads to (\ref{118}).

%%For any integer $0 \leqslant t \leqslant \min\{\,  m, \ n \,\},$  setting
%$r(X) = t,$ $Y = 0$ and $Z=0$ leads to
%r[\phi(X, \,Y, \,Z)] =  r(X) = q;
%$$

(a) for any integer $0  \leqslant  t \leqslant \min\{\, m + q,   \ \  n + p,  \ \  m + n \,\}$
with $m + q \leqslant  \min\{\,  n + p,  \ \  m + n \,\}$ and $m \leqslant p$,   setting
$$
X =0, \ \ \ Y = [\, Y_1, \ 0\,], \ \  Z = [\, Z_1, \ 0\,], \ \ r(Y_1) + r(Z_1) = t
$$
leads to $r[\phi(X, \, Y, \, Z\,)] = r(Y_1) + r(Z_1) = t$;
with $m + q \leqslant  \min\{\,  n + p,  \ \  m + n \,\}$ and $m > p$,
setting
$$
[\,X, \, Y \,] = [\, 0, \ Y_1\,], \ \  Z = [\, Z_1, \ 0\,], \ \ \ r(Y_1) + r(Z_1) = t
$$
leads to $r[\phi(X, \, Y, \, Z\,)] = r(Y_1) + r(Z_1) =  t$;

(b) for any integer $0  \leqslant  t \leqslant \min\{\, m + q,   \ \  n + p,  \ \  m + n \,\}$  with
 $n + p \leqslant  \min\{\,  m + q,  \ \  m + n \,\}$ and $n \leqslant q$, setting
$$
X =0, \ \ \ Y = [\, Y_1, \ 0\,]^T, \ \  Z = [\, Z_1, \ 0\,]^T, \ \ r(Y_1) + r(Z_1) = t
$$
leads to $r[\phi(X, \, Y, \, Z\,)] = r(Y_1) + r(Z_1) = t$; with $n + p \leqslant  \min\{\,  m + q,  \ \  m + n \,\}$
and $n > q$, setting
$$
\left[\!\begin{array}{cc} X \\ Z \end{array}\!\right]  = \left[\!\begin{array}{cc} 0 \\ Z_1
\end{array}\!\right], \ \  Y = \left[\!\begin{array}{cc} Y_1  \\ 0
\end{array}\!\right]
$$
leads to $r[\phi(X, \, Y, \, Z\,)] = r(Y_1) + r(Z_1) = t$;

(c) for any integer $0  \leqslant  t \leqslant \min\{\, m + q,   \ \  n + p,  \ \  m + n \,\}$
with $m + n \leqslant  \min\{\,  m + q,  \ \  n + p \,\}$,   setting
$$
X =0, \ \ \ Y = [\, Y_1, \ 0\,], \ \  Z = [\, Z_1, \ 0\,]^T, \ \  r(Y_1) + r(Z_1) = t
$$
leads to $r[\phi(X, \, Y, \, Z\,)] = r(Y_1) + r(Z_1) = t$, establishing  (\ref{119}). \qquad $\Box$

\begin{lemma} \label{T15}
Let $A \in {\mathbb C}^{m \times n}$ be given$,$ $Y \in \mathbb C^{m \times p},$
$Z \in {\mathbb C}^{q \times n}$  and $U \in {\mathbb C}^{q \times p}$  be three variable matrices$,$ and define
\begin{align}
\phi(Y, \, Z, \, U\,)= \left[\!\begin{array}{cc} A & Y \\ Z & U
\end{array}\!\right]\!.
\label{120}
\end{align}
Then$,$
\begin{align}
\max_{Y \in \mathbb C^{m \times p}, \ Z \in {\mathbb C}^{q \times n}, \, U \in {\mathbb C}^{q \times p}}
r[\phi(X, \,Y, \, U)] & = \min\{\, m + p,  \ \ n + q,  \ \  p + q - r(A)\,\},
\label{121}
\\
\min_{Y \in \mathbb C^{m \times p}, \ Z \in {\mathbb C}^{q \times n}, \, U \in {\mathbb C}^{q \times p}}
r[\phi(X, \,Y, \, U)]  & = r(A).
\label{122}
\end{align}
In particular$,$ for any integer $t$ with $r(A) \leqslant t \leqslant  \min\{\, m + p,  \ \ n + q,  \ \ r(A) + p + q \,\},$
 there exist $Y \in \mathbb C^{m \times p},$ $Z \in {\mathbb C}^{q \times n}$  and $U \in {\mathbb C}^{q \times p}$
 such that
\begin{align}
r[\phi(X, \,Y, \, U)] = t.
\label{123}
\end{align}
\end{lemma}

\noindent {\bf Proof.} \ Without lost generality, we assume that $A$ is given by
\begin{equation}
A = {\rm diag}(\, I_{d},  \ 0 \,).
\label{124}
\end{equation}
Correspondingly,
\begin{align}
\phi(Y,\,Z, \, U) = \left[\!\!\begin{array}{cccccc} I_d & 0 & \widehat{Y}_1
\\
0 & 0 & \widehat{Y}_2
\\
\widehat{Z}_1 & \widehat{Z}_2  & U
\end{array}\!\!\right],
\label{125}
\end{align}
and
\begin{align}
r[\phi(Y,\,Z, \, U)]  = d + r\!\left[\!\!\begin{array}{cc} 0 & \widehat{Y}_2
\\
\widehat{Z}_2 & U - \widehat{Z}_1\widehat{Y}_1 \end{array}\!\!\right]\!.
\label{126}
\end{align}
Applying  (\ref{117}) and (\ref{118}) to the block matrix in (\ref{126}) leads to
\begin{align}
& \max_{Y \in \mathbb C^{m \times p}, \, Z \in \mathbb C^{q \times n}, \, U \in {\mathbb C}^{q \times p}}r\!\left[\!\!\begin{array}{cc}  0 & \widehat{Y}_2
\\
\widehat{Z}_2 & U - \widehat{Z}_1\widehat{Y}_1 \end{array}\!\!\right] = \min\{\, m + p  - r(A), \ \
n +q  - r(A), \ \   p + q - 2r(A)\,\},
\label{127}
\\
& \min_{Y \in \mathbb C^{m \times p}, Z \in \mathbb C^{q \times n}, \, U \in {\mathbb C}^{q \times p}}r\!\left[\!\!\begin{array}{cc}  0 & \widehat{Y}_2
\\
\widehat{Z}_2 & U - \widehat{Z}_1\widehat{Y}_1  \end{array}\!\!\right] = 0.
\label{128}
\end{align}
Substituting (\ref{127}) and (\ref{128}) into (\ref{126}) yields (\ref{121}) and (\ref{122}).
Applying (\ref{119}) to (\ref{127}) and (\ref{128}) leads to (\ref{123}). \qquad $\Box$

\section{Rank optimization of $A + X$}
\renewcommand{\theequation}{\thesection.\arabic{equation}}
\setcounter{section}{2} \setcounter{equation}{0}

One of the special cases in (\ref{11}) is the ordinary sum $A + X$. In this section,
we derive explicit formulas for calculating the extremal ranks of $A + X$ subject to $X$ with a fixed rank.
The formulas obtained will be used in Sections 3.

\begin{theorem}\label{T21}
Let  $A \in {\mathbb C}^{m \times n}$ be given$,$  $X\in {\mathbb C}^{m \times n}$
be a variable matrix$,$  and assume that $s$ and $t$ are two integers satisfying
\begin{align}
0\leqslant s \leqslant t \leqslant \min\{\,m, \ n\,\}.
\label{21}
\end{align}
Then$,$
\begin{enumerate}
\item[{\rm(a)}] The following equalities hold
\begin{align}
\max_{X \in {\mathbb C}_t^{m \times n}}r(\, A + X \,) & = \min\{\, m, \ \  n,  \ \ r(A)  + t\,\},
 \label{22}
\\
\min_{X \in {\mathbb C}_t^{m \times n}}r(\, A + X \,) & = |\, r(A) - t \,|.
\label{23}
\end{align}

\item[{\rm(b)}] The following equalities hold
\begin{align}
\max_{X \in {\mathbb C}^{m \times n},\,s \leqslant r(X) \leqslant t}\!\!\! r(\, A + X \,) & = \min\{\, m, \ \  n, \
r(A)  + t \,\},
 \label{24}
\\
\min_{X \in{\mathbb C}^{m \times n},\, s \leqslant r(X) \leqslant t}\!\!\!r(\, A + X \,) &  = \max\{\, 0, \ s  - r(A), \ r(A) - t\}.
\label{25}
\end{align}

\item[{\rm(c)}] The following equalities hold
\begin{align}
\max_{X \in {\mathbb C}^{m \times n},\,0 \leqslant r(X) \leqslant t} r(\, A + X \,) & = \min\{\, m, \ n, \ r(A)  + t \,\},
 \label{26}
\\
\min_{X \in {\mathbb C}^{m \times n},\,0 \leqslant r(X) \leqslant t}r(\, A + X \,) & = \max\{ \, 0,  \ \  \, r(A) - t  \}.
\label{27}
\end{align}

\item[{\rm(d)}] The following equalities hold
\begin{align}
\max_{X \in {\mathbb C}^{m \times n},\, s \leqslant r(X)\leqslant \min \{ m, \, n\}} r(\, A + X \,) & = \min \{ m, \ \ n\},
 \label{28}
\\
\min_{X \in {\mathbb C}^{m \times n},\, s\leqslant r(X)\leqslant \min \{ m, \, n\}}r(\, A + X \,) & = \max\{ \, 0, \ s  - r(A)  \}.
\label{29}
\end{align}
\end{enumerate}
The matrices $X$ satisfying these equalities can be
formulated from the canonical form of $A.$
\end{theorem}

\noindent {\bf Proof.} \
It is obvious that the right-hand sides of (\ref{22}) and
(\ref{23}) are upper and lower bounds. Without loss of generality, we assume that $A$ is of the form
\begin{align}
A = {\rm diag}(\,I_d,  \ 0  \,).
\label{210}
\end{align}
Let $X = \left[\!\begin{array}{cc} 0 & 0 \\ 0 & I_t \end{array}\!\right]$. If $m  \leqslant
\min\{\, n, \ \ r(A)  + t \,\}$, then $r(\, A + X \,) = m$; if $n  \leqslant
\min \{\, m, \ \ r(A)  + t \,\}$, then $r(\, A + X \,) = n$; if $ r(A) + t \leqslant \min\{\, m, \ \ n \,\}$,
then $r(\, A+ X \,) = r(A) + t,$ so that (\ref{22}) holds.

If $r(A) \leqslant t$, then setting $X =\left[\!\begin{array}{ccc} -I_d & 0 & 0 \\ 0 & I_{t-d} & 0
\\ 0 & 0 & 0 \end{array}\!\right]$ gives $r(\, A + X \,) = t - d$; if $r(A) > t$, then setting $X
=\left[\!\begin{array}{ccc} -I_t & 0 \\ 0 & 0  \end{array}\!\right]$
gives $r(\, A + X \,) = d - t$,  so that (\ref{23}) holds.

Note that
$$
 \{X \in {\mathbb C}^{m \times n} \, | \, s \leqslant r(X) \leqslant t\}
  = {\mathbb C}_s^{m \times n} \cup  {\mathbb C}_{s+1}^{m \times n} \cup \ldots
  \cup {\mathbb C}_t^{m \times n}.
$$
So that
\begin{align}
& \max_{X \in {\mathbb C}^{m \times n},\,s \leqslant r(X) \leqslant t}\!\!\!\!r(\, A + X \,)
 = \max\left\{ \max_{X\in {\mathbb C}_s^{m \times n}}\!\!\!r(\, A + X \,), \
 \max_{X\in {\mathbb C}_{s+1}^{m \times n}}\!\!\!r(\, A + X \,), \ldots,
  \max_{X\in {\mathbb C}_t^{m \times n}} \!\!\!r(\, A + X \,) \right\},
 \label{211}
\\
 & \min_{X \in {\mathbb C}^{m \times n},\,s \leqslant r(X) \leqslant t}\!\!\!\!r(\, A + X \,)
  = \min \left\{ \min_{X\in{\mathbb C}_s^{m \times n}}\!\!r(\, A + X \,), \
\min_{X\in {\mathbb C}_{s+1}^{m \times n}}\!\!r(\, A + X \,),  \ldots,
\min_{X\in {\mathbb C}_{t}^{m \times n}} \!\!r(\, A + X \,)\right\}.
\label{212}
\end{align}
Substituting (\ref{22}) and (\ref{23}) for  $r(X) = s,\,  s +1, \ldots, \, t$ into (\ref{211}) and (\ref{212}) and
making the max-min comparison, we obtain
\begin{align}
& \max_{X \in {\mathbb C}^{m \times n},\,s \leqslant r(X) \leqslant t} \!\!r(\, A + X \,)  \nb
\\
& = \max\left\{ \min\{\, m, \ n, \
r(A)  + s\,\}, \ \min\{\, m, \ n, \  r(A)  + s +1\,\}, \, \ldots, \,  \min\{\, m, \ n, \ r(A)  + t\,\}
\right\} \nb
\\
& = \min\{\, m, \ n, \  r(A)  + t\,\},
 \label{213}
\\
& \min_{X \in {\mathbb C}^{m \times n},\, s \leqslant r(X) \leqslant t}\!\!r(\, A + X \,)  \nb
\\
& = \min \left\{ |\, r(A) - s \,|,  \ |\, r(A) - s -1 \,|, \, \ldots, \, |\, r(A) - t \,| \right\}
 = \max\{\, 0, \ s  - r(A), \ r(A) - t  \},
\label{214}
\end{align}
establishing (\ref{24}) and (\ref{25}), as well as (\ref{26})--(\ref{29}). \qquad $\Box$

\medskip

The results in the section show that the matrix rank optimization formulated in (\ref{12})--(\ref{15}) are combinatorial in nature.

\section{Rank optimization of $A+ BXC$}
\renewcommand{\theequation}{\thesection.\arabic{equation}}
\setcounter{section}{3} \setcounter{equation}{0}

A matrix-valued function for complex matrices is a map between matrix spaces ${\mathbb C}^{m \times n}$ and
${\mathbb C}^{p \times q}$, which can generally be written as
$$
Y = f(X) \ \  {\rm for}  \ \ Y \in {\mathbb C}^{m \times n} \ \ {\rm and} \ \  X \in {\mathbb C}^{p \times q},
$$
 or briefly, $f:  {\mathbb C}^{m \times n} \rightarrow {\mathbb C}^{p \times q}$. Eq. (\ref{11}) is in fact the simplest case
 of all matrix-valued functions, which is extensively studied from theoretical and applied points of view.

According to \cite{Zha-1}, substituting  (\ref{112}) into (\ref{11}) yields
\begin{equation}
 \phi(X)=P\Sigma_AQ +  P\Sigma_BUXV\Sigma_CQ  =P(\,\Sigma_A  + \Sigma_BUXV\Sigma_C\,)Q,
\label{31}
\end{equation}
which we call a canonical form of (\ref{11}).  Many properties of the matrix-valued function $A + BXC$  can be
derived from the canonical form. For instance, the rank of $A + BXC$ is
\begin{equation}
 r(\, A + BXC \,)=r(\,\Sigma_A  + \Sigma_BY\Sigma_C\,), \ \ r(X) = r(Y),
 \label{32}
\end{equation}
where $Y=UXV$. Partition it as
\begin{align}
& \begin{array}{c@{\hspace{4pt}}l}
Y =\left[\!\! \begin{array}{cccc}
\qquad  Y_{11} \quad  & Y_{12}  & Y_{13} & Y_{14}  \\
\qquad  Y_{21} \quad  & Y_{22}  & Y_{23} & Y_{24}  \\
\qquad  Y_{31} \quad   & Y_{32} & Y_{33} & Y_{34}  \\
\qquad  Y_{41} \quad   & Y_{42} & Y_{43}  & Y_{44}
 \end{array}
\!\!\right] &
\begin{array}l
j \\  p - j -u - s_2  \\ u \\ s_2
\end{array}
\\
{\footnotesize  \quad   q- l - u - s_1 \quad  l  \qquad \  u \qquad  \ s_1}
\end{array}\!\!\!.
\label{33}
\end{align}
Then we have
\begin{align}
\Sigma_A  + \Sigma_BY\Sigma_C = \left[\!\!
\begin{array}{cccccc}
 I_j & 0 & Y_{12}  &  Y_{13}S_C  & Y_{14} & 0
 \\
 0    &     I_k    & 0     & 0  & 0  & 0
 \\
 0  &  0  &  I_l  & 0 & 0  & 0
 \\
 0  &   0  &  S_BY_{32}  & S_A + S_BY_{33}S_C  & S_BY_{34} & 0
 \\
  0 & 0  & Y_{42}  &  Y_{43}S_C  & Y_{44}  & 0
 \\
 0 &  0  & 0 & 0  & 0 & 0
 \end{array}
\!\!\right],
\label{34}
\end{align}
and
\begin{align}
r(\, \Sigma_A  + \Sigma_BY\Sigma_C \,)&  =  j + k + l + r\!\left[\!\begin{array}{cc} S_A + S_BY_{33}S_C  & S_BY_{34}
\\
 Y_{43}S_C  & Y_{44}  \end{array}\!\right]  \nb
\\
&= r[\, A, \, B\,]+ r\!\left[\!\begin{array}{c} A \\ C
\end{array}\!\right] - r\!\left[\!\begin{array}{cc} A & B \\ C& 0
\end{array}\!\right] + r\!\left[\!\begin{array}{cc}
S_B^{-1}S_AS_C^{-1} +  Y_{33} & Y_{34}\\ Y_{43} & Y_{44} \end{array}\!\right]\!.
\label{35}
\end{align}
So that
\begin{align}
&  \max_{Y \in {\mathbb C}_t^{p \times q}}r(\, \Sigma_A  + \Sigma_BY\Sigma_C \,)
 = r[\, A, \, B\,]+ r\!\left[\!\begin{array}{c} A \\ C
\end{array}\!\right] - r\!\left[\!\begin{array}{cc} A & B \\ C& 0
\end{array}\!\right] + \max_{Y \in {\mathbb C}_t^{p \times q}}r(\, S +  \widehat{Y}\,),
\label{36}
\\
& \min_{Y\in {\mathbb C}_t^{p \times q}}r(\, \Sigma_A  + \Sigma_BY\Sigma_C \,)
 = r[\, A, \, B\,]+ r\!\left[\!\begin{array}{c} A \\ C
\end{array}\!\right] - r\!\left[\!\begin{array}{cc} A & B \\ C& 0
\end{array}\!\right] + \min_{Y \in {\mathbb C}_t^{p \times q}}r(\, S +  \widehat{Y}\,),
\label{37}
\end{align}
where
$$
S = \left[\!\begin{array}{cc}  S_B^{-1}S_AS_C^{-1} &  0 \\ 0 & 0 \end{array}\!\right]\!, \ \
\widehat{Y} = \left[\!\begin{array}{cc} Y_{33} & Y_{34}\\ Y_{43} & Y_{44} \end{array}\!\right]\!.
$$
Applying Theorem \ref{21}(b) to (\ref{36}) and (\ref{37}) yields the main result in the note.

\begin{theorem} \label{T31}
Let $\phi(X)$ be as given in {\rm (\ref{11})}$,$  and assume that  $t$ is an integer satisfying
$0\leqslant t \leqslant \min\{\,p, \ \ q\,\}.$ Also define
\begin{align}
G =[\,A, \, B\,], \ \  H = \!\left[\!\begin{array}{c} A \\ C
 \end{array}\!\right]\!, \ \  M = \left[\!\begin{array}{cc} A  & B\\ C & 0 \end{array}\!\right]\!.
\label{38}
\end{align}
Then$,$
\begin{align}
&  \max_{X \in {\mathbb C}_t^{p \times q}}\!\!r(\, A + BXC \,)
 =\min\left\{ r[\,A, \, B\,] \ \ r\!\left[\!\begin{array}{c} A \\ C
 \end{array}\!\right]\!, \ \ r(A) + t \right\},
\label{39}
\\
 & \min_{ X \in {\mathbb C}_t^{p \times q}}\!\!r(\, A + BXC \,)
 = \max \left\{ r[\,A, \, B\,] + r\!\left[\!\begin{array}{c} A \\ C
 \end{array}\!\right] - r\!\left[\!\begin{array}{cc} A  & B\\ C & 0 \end{array}\!\right]\!,
 \ \ r[\,A, \, B\,] + r\!\left[\!\begin{array}{c} A \\ C
 \end{array}\!\right] - r(A) + t - p - q, \ \ r(A) - t \right\}.
\label{310}
\end{align}
In consequences$,$
\begin{enumerate}
\item[{\rm(a)}] Under $m = n,$ there exists an $X\in {\mathbb C}_t^{p \times q}$
such that $A + BXC$ is nonsingular if  and only if
\begin{align}
r(G) = r(H) = m \ \  and  \ \ r(A) \geqslant m -t.
\label{311}
\end{align}

\item[{\rm(b)}] There exists an $X\in {\mathbb C}_t^{p \times q}$
such that $A + BXC = 0$ if and only if
\begin{align}
{\mathscr R}(A) \subseteq {\mathscr R}(B), \ \ {\mathscr R}(A^*) \subseteq {\mathscr R}(C^*), \ \
r(G) + r(H) \leqslant r(A) - t + p + q, \ \  r(A) \leqslant t.
\label{312}
\end{align}

\item[{\rm(c)}] Under $t \neq 0,$ the rank of $A + BXC$ is invariant for all $X\in {\mathbb C}_t^{p \times q}$
 if and only if one of the following conditions holds$:$
  \begin{enumerate}
\item[{\rm(i)}] $r(M) = r(G),$

 \item[{\rm(ii)}] $r(M) = r(H),$

 \item[{\rm(iii)}] $r(M) = r(G) + f(H) - r(A) - t,$

 \item[{\rm(iv)}]  $r(G) = r(A) + p + q - t$,

\item[{\rm(v)}]  $r(H)  = r(A) +  p + q - t$,

\item[{\rm(vi)}]  $r(G) +r(H) = 2r(A) + p + q,$ namely$,$ $r(B) = p,$ $r(C) = q,$ ${\mathscr R}(A)\cap {\mathscr R}(B) = \{0\}$
and ${\mathscr R}(A^*)\cap {\mathscr R}(C^*) = \{0\}.$
\end{enumerate}

\item[{\rm(d)}] Under ${\mathscr R}(A) \subseteq {\mathscr R}(B)$ and ${\mathscr R}(A^*) \subseteq {\mathscr R}(C^*),$
\begin{align}
&  \max_{X \in {\mathbb C}_t^{p \times q}}\!\!r(\, A + BXC \,)
 =\min\left\{ r(B), \ \ r(C), \ \ r(A) + t \right\},
\label{39a}
\\
 & \min_{ X \in {\mathbb C}_t^{p \times q}}\!\!r(\, A + BXC \,)
 = \max \left\{ 0, \ \ r(B) + r(C) - r(A) + t - p - q, \ r(A) - t \right\}.
\label{310a}
\end{align}
\end{enumerate}
\end{theorem}

\noindent {\bf Proof.} \ Let $z = p + q - 2r(M) + r(G) +  r(H).$
Then we  find by (\ref{123}),  (\ref{24}) and (\ref{25}) that
\begin{align}
& \max_{Y \in {\mathbb C}_t^{p \times q}}r(\, S +  \widehat{Y}\,)  =
\max_{t - z \leqslant  r(\widehat{Y})\leqslant t }r(\, S +  \widehat{Y}\,)
= \min\{\, u + s_1, \ \ u + s_2, \ \ r(S) + t \,\} \nb
\\
& = \min\left\{\, r(M) -  r(G), \ \ r(M) - r(H) , \ \ r(M) + r(A) -  r(G) -  r(H) + t \right\},
\label{313}
\\
& \min_{Y \in {\mathbb C}_t^{p \times q}}r(\, S +  \widehat{Y}\,)  =
\min_{t - z \leqslant  r(\widehat{Y})\leqslant t }r(\, S +  \widehat{Y}\,)  =
\max\{\, 0, \ \  t - z - r(S), \ \ r(S)  - t \,\} \nb
\\
& = \max\left\{\, 0, \ \  t - p - q + r(M) - r(A),  \ \ r(M) + r(A) -  r(G) -  r(H) - t  \right\}.
\label{314}
\end{align}
Substituting (\ref{313}) and (\ref{314}) into (\ref{36}) and (\ref{37}) yields
(\ref{39}) and (\ref{310}).  Setting (\ref{39}) equal to $m$ yields (\ref{311});
setting (\ref{310}) equal to $0$ yields (\ref{312});  setting (\ref{39}) equal to
(\ref{310}) yields the results  in (c). \qquad  $\Box$

\medskip

Eqs. (\ref{39}) and (\ref{310}) show that the extremal ranks of
 (\ref{11}) can be calculated  exactly  from the two formulas without
 knowing how to choose the feasible matrices $X$. So that they can be used
 independently in describing behaviors of $A + BXC$, as shown in Theorem \ref{T31}(a)--(d).

Recall that any matrix $X\in {\mathbb C}_t^{p \times q}$   can be written as
 a product $X = YZ$, where $\in {\mathbb C}_t^{p \times t}$ and $Z\in {\mathbb C}_t^{t \times q}$
 with  $r(Y) =  r(Z) = t$. So that (\ref{39}) and (\ref{310})  can be represented as follows.

\begin{corollary} \label{T32}
Let $\phi(X)$ be as given in {\rm (\ref{11})}$,$ $t$ be an integer satisfying
$0\leqslant t \leqslant \min\{\,p, \ \ q\,\},$ and $G,$ $H$ and $M$ be the matrices in
 {\rm (\ref{38})}$.$  Then$,$
\begin{align}
 \max_{Y \in {\mathbb C}_t^{p \times t},\, Z \in {\mathbb C}_t^{t \times q}}\!\!\!\!r(\, A + BYZC \,)
 &=\min\left\{ r(G), \ \ r(H), \ \ r(A) + t \right\},
\label{ww317}
\\
 \min_{Y \in {\mathbb C}_t^{p \times t},\, Z \in {\mathbb C}_t^{t \times q}}\!\!\!\!r(\, A + BYZC \,)
 &= \max \left\{ r(G)+ r(H) - r(M), \ \ r(G) + r(H) - r(A) + t - p - q,  \ \ r(A) - t \right\}.
\label{ww318}
\end{align}
\end{corollary}

\begin{corollary} \label{T33}
Let $\phi(X)$ be as given in {\rm (\ref{11})}$,$ $G,$ $H$ and $M$ be the matrices in
 {\rm (\ref{38})}$,$ and assume that $s$ and $t$ are two integers satisfying
\begin{align}
0\leqslant s \leqslant t \leqslant \min\{\, p, \ \ q\,\}.
\label{315}
\end{align}
Then$,$
\begin{align}
&  \max_{X \in \mathbb C^{p \times q},  \, s \leqslant r(X) \leqslant t}\!\!r(\, A + BXC \,)
 =\min\left\{ r(G), \ \ r(H), \ \ r(A) + t \right\},
\label{316}
\\
 & \min_{ X \in \mathbb C^{p \times q},  \, s \leqslant r(X) \leqslant t}\!\!r(\, A + BXC \,)
 = \min\{u_s, \ u_{s+1}, \ldots, u_{t} \},
 \label{317}
\end{align}
where
$$
u_l= \max \left\{ r(G)+ r(H) - r(M), \ \ r(G) + r(H) - r(A) + l - p - q, \ \ r(A) - l \right\}, \ \
 l = s, \ s+1, \ldots, t.
$$
\end{corollary}

\begin{corollary} \label{T34}
Let $A \in {\mathbb C}^{m \times n},$ $B \in {\mathbb C}^{m \times p}$
 and $C \in {\mathbb C}^{p \times n}$ be given$.$ Then$,$
\begin{align}
&  \max_{X \in {\mathbb C}_p^{p \times p}}\!\!r(\, A + BXC \,)
 =\min\left\{ r[\, A, \, B\,], \ \  r\!\left[\!\begin{array}{c} A \\ C
\end{array}\!\right] \right\},
\label{318}
\\
 & \min_{ X \in {\mathbb C}_p^{p \times p}}\!\!r(\, A + BXC \,)
 = \max \left\{ r[\, A, \, B\,]+ r\!\left[\!\begin{array}{c} A \\ C
\end{array}\!\right] - r\!\left[\!\begin{array}{cc} A & B \\ C& 0
\end{array}\!\right], \ \ r[\, A, \, B\,]+ r\!\left[\!\begin{array}{c} A \\ C
\end{array}\!\right]  - r(A) - p \right\}.
\label{319}
\end{align}
\end{corollary}

\begin{corollary} \label{T35}
Let $0 \neq B \in {\mathbb C}^{m \times p}$ and $0 \neq C \in {\mathbb C}^{q \times n}$ be given$,$
and assume that  $t$ is an integer satisfying $1\leqslant t \leqslant \min\{\,p, \ \ q\,\}.$
 Then$,$
\begin{align}
&  \max_{X \in \mathbb C^{p \times q}, \,r(X) = t}\!\!r(\,BXC \,)
 =\min\left\{ r(B), \ \ r(C), \ \  t \right\},
\label{320}
\\
 & \min_{ X \in \mathbb C^{p \times q},  \,r(X) = t}\!\!r(\,BXC \,)
 = \max \left\{ 0, \ \ r(B) + r(C)  + t - p - q \right\}.
\label{321}
\end{align}
In consequences$,$
\begin{enumerate}
\item[{\rm(i)}] Under $m = n,$ there exists an $X\in {\mathbb C}^{p \times q}$ with
$r(X) =t$ such that $BXC$ is nonsingular if  and only if
\begin{align}
r(B) = r(C) = m \ \  and  \ \ t \geqslant m.
\label{322}
\end{align}

\item[{\rm(ii)}] There exists an $X\in {\mathbb C}^{p \times q}$ with $r(X)= t$
such that $BXC = 0$ if and only if
\begin{align}
r(B) + r(C) \leqslant  p + q - t.
\label{323}
\end{align}

\item[{\rm(iii)}] The rank of $BXC$ is invariant for all $X\in {\mathbb C}^{p \times q}$ with $r(X)= t$
if and only if
\begin{align}
r(B) =  p + q - t, \ \ or \ \ r(C) =  p + q - t,  \ \ or \ \ r(B) =  p  \ and  \  r(C) =  q.
\label{324}
\end{align}
\end{enumerate}
\end{corollary}

\section{Completing a partially-specified block matrix with extremal ranks}
\renewcommand{\theequation}{\thesection.\arabic{equation}}
\setcounter{section}{4} \setcounter{equation}{0}

As an application of the results in the previous section, we consider the rank of the
following partially specified block matrix
\begin{align}
\phi(X) =\left[\!\begin{array}{cc} A & B \\ C & X
\end{array}\!\right]  =  \left[\!\begin{array}{cc} A  & B \\ C & 0 \end{array}\!\right]  +
\left[\!\begin{array}{c} 0\\ I_q
\end{array}\!\right]X[\, 0, \ I_p\,],
\label{41}
\end{align}
where $A \in  {\mathbb C}^{m \times n}$, $B \in  {\mathbb C}^{m \times p}$ and $C\in \mathbb C^{q
\times n}$ are given, and $X\in  {\mathbb C}^{q \times p}$ is a variable matrix, which obviously is
a special case of (\ref{11}).  Conversely, the rank of (\ref{11}) can equivalently be written as
$$
r(\, A + BXC \,) =  r\!\left[\!\begin{array}{ccc} A & B  &  0\\ C & 0 & I_q \\
0 & I_p & X \end{array}\!\right]  - p - q,
$$
the block matrix in which is a special case of (\ref{41}) as well.

\begin{theorem} \label{T41}
Let $\phi(X)$ be as given in {\rm (\ref{41})}$,$  and assume that $s$ and $t$ are two integers satisfying
\begin{align}
0\leqslant s \leqslant t \leqslant \min\{\,p, \ q\,\}.
\label{42}
\end{align}
Also define
\begin{align}
G =[\,A, \, B\,], \ \  H = \!\left[\!\begin{array}{c} A \\ C
 \end{array}\!\right], \ \  M = \left[\!\begin{array}{cc} A  & B\\ C & 0 \end{array}\!\right].
\label{43}
\end{align}
Then$,$
\begin{enumerate}
\item[{\rm(a)}] The following equalities hold
\begin{align}
\max_{X \in {\mathbb C}_t^{q \times p}}\!\!r[\phi(X)]  &= \min\left\{  r(G) + q, \ \ r(H) + p, \ \ r(M) + t \right\},
 \label{44}
\\
\min_{X\in {\mathbb C}_t^{q \times p}}\!\!r[\phi(X)] &= \max\left\{ r(G) + r(H) - r(A), \ \ r(G) + r(H) -r(M) + t, \ \ r(M) - t \right\}.
 \label{45}
\end{align}
In consequence$,$
\begin{enumerate}
\item[{\rm(i)}] Under $m + q = n +p,$ there exists an $X \in {\mathbb C}_t^{q \times p}$ such that $\phi(X)$ is nonsingular
if and only if $r(G) =m,$ $r(H) =n$ and $r(M)\geqslant m + q -t.$

\item[{\rm(ii)}] Under $m + q = n +p,$  $\phi(X)$ is nonsingular for all $X \in {\mathbb C}_t^{q \times p}$
if and only if $r(G) + r(H) - r(A) =m +q$ or $r(G) + r(H) -r(M) = m +q - t,$ or $r(M) = m +q$ and $t =0.$
\end{enumerate}

\item[{\rm(b)}] Under ${\mathscr R}(B) \subseteq {\mathscr R}(A)$ and ${\mathscr R}(C^*) \subseteq {\mathscr R}(A^*),$ the following
equalities hold
\begin{align}
& \max_{X \in {\mathbb C}_t^{q \times p}}\!\!r[\phi(X)] = \min\left\{  r(A) + q, \ \ r(A) + p, \ \ r(M) + t \right\},
 \label{46}
\\
& \min_{X\in {\mathbb C}_t^{q \times p}} \!\!r[\phi(X)] = \max\left\{ r(A), \ \ 2r(A) -r(M) + t, \ \ r(M) - t \right\}.
 \label{47}
\end{align}

\item[{\rm(c)}] Under ${\mathscr R}(A) \subseteq {\mathscr R}(B)$ and ${\mathscr R}(A^*) \subseteq {\mathscr R}(C^*),$
the following equalities hold
\begin{align}
\max_{X \in {\mathbb C}_t^{q \times p}}\!\!r[\phi(X)]  &= \min\left\{  r(B) + q, \ \ r(C) + p, \ \ r(B) + r(C)+ t \right\},
 \label{48}
\\
\min_{X\in {\mathbb C}_t^{q \times p}}\!\!r[\phi(X)] &= \max\left\{ r(B) + r(C) - r(A), \ \  t, \ \ r(B) + r(C) - t \right\}.
 \label{49}
\end{align}

\item[{\rm(d)}] The following equalities hold
\begin{align}
\max_{X \in {\mathbb C}^{q \times p},\, s\leqslant r(X) \leqslant t}\!\!r[\phi(X)]  &= \min\left\{  r(G) + q, \ \ r(H) + p, \ \ r(M) + t \right\},
 \label{410}
\\
\min_{X\in {\mathbb C}^{q \times p},\,s\leqslant r(X) \leqslant t}\!\!r[\phi(X)] &= \min\{u_s, \ u_{s+1}, \ldots, u_{t} \},
 \label{411}
\end{align}
where
$$
u_l= \max \left\{ r(G)+ r(H) - r(A), \ \ r(G) + r(H) - r(M) + l, \ \ r(A) - l \right\}, \ \
 l = s, \ s+1, \ldots, t.
$$
\end{enumerate}
\end{theorem}

\begin{theorem} \label{T42}
Let \begin{align}
\phi(X) =\left[\!\begin{array}{cc} A - X & B - X  \\ C - X &  D - X
\end{array}\!\right] = \left[\!\begin{array}{cc} A  & B \\ C & D \end{array}\!\right]  -
\left[\!\begin{array}{c} I_m \\ I_m
\end{array}\!\right]X[\, I_n, \ I_n\,],
\label{412}
\end{align}
where $A, \, B, \, C, \, D \in \mathbb C^{m \times n}.$ Also define
\begin{align}
G =[\,A - C, \, B - D\,], \ \  H = \!\left[\!\begin{array}{c} A - B \\ C - D
 \end{array}\!\right], \ \  M = \left[\!\begin{array}{cc} A  & B\\ C & D \end{array}\!\right].
\label{413}
\end{align}
Then the following equalities hold
\begin{align}
\max_{X \in {\mathbb C}_t^{m \times n}}\!\!r[\phi(X)]  &= \min\left\{  r(G) + m, \ \ r(H) + n, \ \ r(M) + t \right\},
 \label{414}
\\
\min_{X\in {\mathbb C}_t^{m \times n}}\!\!r[\phi(X)] &= \max\left\{ r(G) + r(H) - r(A - B - C + D), \ \ r(G) + r(H) -r(M) + t, \ \ r(M) - t \right\}.
 \label{415}
\end{align}
\end{theorem}

Because the right-hand sides of (\ref{44}), (\ref{45}),  (\ref{414}) and (\ref{415}) can be calculated exactly, these results
can be taken as test examples in fixed-rank or lower-rank approximation and perturbation analysis
of matrices. They can also be used to verify the correctness the effectiveness of various numerical algorithms in
rank minimization problems occurred in recent years.

\section{Concluding remarks}
\renewcommand{\theequation}{\thesection.\arabic{equation}}
\setcounter{section}{5} \setcounter{equation}{0}

 Closed-form formulas are established for calculating the extremal ranks in (\ref{12})--(\ref{15}).
This work shows a surprising fact that many matrix rank optimization problems do exist
analytical solutions for calculating the extremal ranks.

Besides (\ref{12})--(\ref{15}), a more popular problem  is to  minimize the norm of (\ref{11}) subject to low-rank  constraint
\begin{align}
& \text{minimize} \ \|\, A + BXC \,\|_F  \ \ \ \ \ \ \ \  \mbox{s.t.}  \ \  X \in  {\mathbb C}^{p\times q}
\ \ \mbox{and} \ \   r(X) \leqslant t,
\label{13a}
\end{align}
see \cite{FT,Son}. So that a comparison of solutions to (\ref{11})  and  (\ref{13a}) can further be discussed.

%The ultimate goal of optimization theory is to seek analytical solutions of given optimization problems. %Once analytical solutions are obtained, people can use them to develop various

Matrix rank optimization problem is really a fruitful research field in both matrix analysis and optimization theory.
In recent years, many numerical methods were developed for matrix rank minimization problems based on approximation
and iteration methods. However, there is no  evidence that these numerical methods can make the matrix-rank-objective
functions really attain their minimal values. Because the exact extremal ranks of $A + BXC$ can be calculated
by the analytical formulas in this note, they will set a principle for verifying the correctness and effectiveness
of these numerical methods.\\


\begin{thebibliography}{99}
\def\itemsep{-0.5ex}
\small

\bibitem{DG}  B. De Moor and G.H. Golub, The restricted singular value decomposition: properties
and applications. SIAM J. Matrix Anal. Appl. 12(1991), 401--425.

\bibitem{EY} C. Eckart and G. Young, The approximation of one matrix by another of lower rank. Psychometrika 19(1936), 211--218.

\bibitem{FT}  S. Friedland and A. Torokhti, Generalized rank-
constrained matrix approximations. SIAM J. Matrix Anal. Appl. 29(2007), 656--659.

\bibitem{LT-lama} Y. Liu and Y. Tian,  How to use RSVD to solve the matrix equation $A=BXC$.
  Linear Multilinear Algebra 58(2010), 537--543.



\bibitem{MGC} S. Ma, D. Goldfarb, L. Chen, Fixed point and Bregman iterative methods for matrix rank minimization,
Math. Program. Ser. A 128(2011), 321--353.


\bibitem{Sc}  E. Schmidt,  Zur Theorie der linearen und nonlinearen Integralgleichungen. Math. Ann. 63(1907), 433--476.

\bibitem{Son}  D.  Sondermann, Best approximate solutions to ma-
trix equations under rank restrictions. Statistical Papers 27(1986), 57--66.


\bibitem{Ti-1} Y.~Tian, The maximal and minimal ranks of  some expressions
of generalized inverses of matrices.  Southeast Asian Bull. Math.
25(2002), 745--755.

\bibitem{TC}  Y.~Tian and S.~Cheng,  The maximal and
minimal ranks of $A - BXC$ with applications.  New York J. Math.
9(2003), 345--362.

\bibitem{Zha-1} H.~Zha, The restricted singular value decomposition of matrix triplets.
SIAM. J. Matrix. Anal. Appl. 12(1991), 172--194.

\bibitem{Zha-2} H.~Zha, A numerical algorithm for computing the
restricted singular value decomposition of matrix triplets. Linear Algebra Appl. 168(1992), 1--25.
\end{thebibliography}
\end{document}